\documentclass[a4paper,11pt,leqno,twoside]{article}

\usepackage{amsmath, amsthm, amssymb}

\numberwithin{equation}{section}

\def\bZ{\mathbb{Z}}

\def\bR{\mathbb{R}}

\def\bK{\mathbb{K}}

\def\bJ{\mathbb{J}}

\def\bE{\mathbb{E}}

\def\mP{\mathcal{P}}
\def\d{\mathrm{d}}
\def\PP{\mathrm{P}}
\def\mD{\mathcal{D}}
\def\mS{\mathfrak{S}}
\def\be{\mathbf{e}}
\def\bv{\mathbf{v}}
\def\Ad{\mathrm{Ad}}

\def\tpsi{\widetilde{\psi}}
\def\tPsi{\widetilde{\Psi}}

\def\Ai{\mathrm{Ai}}
\def\Airy{\mathrm{Airy}}

\def\Pf{\mathrm{Pf}}
\def\odd{\mathrm{odd}}

\def\Sh{\mathrm{Sh}}

\def\({ \left( }
\def\){ \right)}
\def\[{ \left[ }
\def\]{ \right]}

\theoremstyle{plain}
\newtheorem{thm}{Theorem}[section]
\newtheorem{prop}[thm]{Proposition}
\newtheorem{lem}[thm]{Lemma}
\newtheorem{cor}[thm]{Corollary}
\theoremstyle{definition}

\newtheorem{example}{Example}[section]
\theoremstyle{remark}
\newtheorem{remark}{Remark}[section]
\theoremstyle{conjecture}

\newcommand{\textem}[1]{{\bfseries #1}}

\title{\bfseries Correlation functions of \\ the shifted Schur measure}
\author{\textsc{Sho Matsumoto}}

\date{\today}

\begin{document}
\setlength{\baselineskip}{16pt}
\maketitle

 
%
\begin{abstract}
The shifted Schur measure introduced in \cite{TracyWidom200?} is a measure on the set of all strict partitions
$\lambda=(\lambda_1> \lambda_2 > \dots >\lambda_{\ell} >0)$,
which is defined by Schur $Q$-functions.
The main aim of this paper is to calculate
the correlation function of this measure, which is given by a pfaffian.
As an application, 
we prove that a limit distribution of $\lambda_j$'s with respect to 
a shifted version of the Plancherel measure for symmetric groups 
is identical with the corresponding distribution of the original Plancherel measure 
(\cite{BDJ,BOO,JohanssonDis,Okounkov2000}).
In particular, we obtain a limit distribution of the length of the longest ascent pair for a random permutation.
Further we give expressions of the mean value and the variance of the size $|\lambda|$ with respect to the measure
defined by Hall-Littlewood functions.

\par\noindent
\textem{2000 Mathematics Subject Classification} : Primary 60C05; Secondary 05E05.
\par\noindent
\textem{Key Words and Phrases} : shifted Schur measure, Schur $Q$-functions, correlation functions, limit distributions,
Plancherel measures, ascent pair, random permutations, Hall-Littlewood functions, Tracy-Widom distribution.
\end{abstract}

%
\section{Introduction}
%
Let $\pi$ be a permutation in the symmetric group $\mS_N$ 
and $\ell(\pi)$ the length of the longest increasing subsequence in $\pi$.
Concerning a limit distribution of $\ell(\pi)$ with respect to
the uniform measure $\PP_{\mathrm{uniform}, N}$ on $\mS_N$, 
it is proved in \cite{BDJ} that
\begin{equation} \label{thmBDJ}
\lim_{N \to \infty} \PP_{\mathrm{uniform},N} \( \frac{ \ell(\pi) - 2 \sqrt{N} }{ N^{1/6}} <s \) = F_2(s),
\end{equation}
where $F_2(s)$ is the Tracy-Widom distribution.
The Tracy-Widom distribution is defined by the Fredholm determinant for the Airy kernel.
Namely,
let $\Ai (x)$ be the Airy function
\begin{equation} \label{Airyfunction}
   \Ai (x)= \frac{1}{2 \pi \sqrt{-1}} \int_{\infty e^{-\pi \sqrt{-1}/3}}^{\infty e^{\pi \sqrt{-1}/3}}        
                     e^{z^3/3 -x z} \d z
\end{equation}
and $K_{\mathrm{Airy}}(x,y)$ the Airy kernel
\begin{equation} \label{Airykernel}
   K_{ \mathrm{Airy} } (x,y) =  \int_0^{\infty} \Ai(x+z) \Ai (z+y) \d z.
\end{equation}
Then the Tracy-Widom distribution $F_2(s)$ is defined by 
\begin{align} 
F_2(s) &= \det(I- K_{\Airy}) |_{L^2([s,\infty))} \label{TWdist} \\
&= 1+ \sum_{k=1}^{\infty} \frac{(-1)^k}{k!} \int_{[s,\infty)^k} \det (K_{\Airy} (x_i,x_j))_{1 \leq i,j \leq k}
 \d x_1 \dots \d x_k \notag
\end{align} 
and gives a limit distribution of the scaled largest eigenvalue of 
a Hermitian matrix from the Gaussian Unitary Ensemble (GUE), see \cite{TracyWidom1994}.

As we see below,
the Plancherel measure for partitions is related to the distribution of the length $\ell(\pi)$.
Let $f^{\lambda}$ be the number of standard tableaux of shape $\lambda$.
The Plancherel measure assigns to each $\lambda \vdash N$ the probability
\begin{equation} \label{Plancherel}
\PP_{\mathrm{Plan},N} (\{\lambda \}) = \frac{(f^{\lambda})^2}{N!}.
\end{equation}
Then it follows from the Robinson-Schensted correspondence that
\begin{equation} \label{uniformPlancherel}
\PP_{\mathrm{uniform},N} (\{ \pi \in \mS_N | \ell(\pi) = h \})
= \PP_{\mathrm{Plan},N} ( \{ \lambda \in \mP_N | \lambda_1 = h \} ),
\end{equation}
where $\mP_N$ is the set of all partitions of $N$ (see e.g. \cite{SaganBook}).
Hence the equation \eqref{thmBDJ} also describes a limit distribution of $\lambda_1$ with respect to Plancherel measures.
This result has been extended in \cite{BOO,JohanssonDis,Okounkov2000} 
to the other rows $\lambda_j$'s in a general position of a partition.
The key of the proof in \cite{BOO} is
a calculation of {\it correlation functions} of the poissonization of the Plancherel measures.
We can see the other asymptotics with respect to the Plancherel measure in e.g. \cite{Hora}. 

On the other hand,
the Schur measure introduced in \cite{Okounkov2001} is a measure which assigns
to each partition the product of two Schur functions.
Okounkov \cite{Okounkov2001} calculated the correlation function of the Schur measure by using the infinite wedge.
The correlation function of the poissonized Plancherel measure is obtained as a specialization of the one
for the Schur measure.

The main aim of this paper is to calculate the correlation function of the shifted Schur measure (see Theorem \ref{corSS}).
The shifted Schur measure, introduced in \cite{TracyWidom200?}, is a measure on the set of all strict partitions,
which is defined by Schur $Q$-functions instead of Schur functions.
The correlation function is expressed as a pfaffian and is actually
calculated by operators on the exterior algebra
in place of the infinite wedge in \cite{Okounkov2001}.
Further, as an application, we obtain a shifted version of the corresponding result for a limit distribution
of $\lambda_j$'s
in \cite{BOO,JohanssonDis,Okounkov2000} (see Theorem \ref{mainthm}).
In particular, we find that a limit distribution of the length of the longest ascent pair for a random permutation 
is given by the Tracy-Widom distribution (see Corollary \ref{maincor}).
Since the proof is similar to the one in \cite{BOO},
we only discuss its main point.  

In the final section, we study about a measure defined by Hall-Littlewood functions.
The measure is considered as a natural extension of the Schur measure and the shifted Schur measure.
We obtain expressions of the mean value $\bE(|\lambda|)$ 
and the variance $\mathrm{Var} (|\lambda|)$ of the size $|\lambda|$
with respect to this measure explicitly.
Actually, each value is written as a sum of the product of certain power-sum functions (see Theorem \ref{MeanValue}). 
This expression of $\bE(|\lambda|)$ naturally leads us a similar study of $\bE(\lambda_1)$. 
By observing various examples,
in the end of the section, 
we remark that there is a certain common property of expressions of $\bE(\lambda_1)$ among these examples.

%
\section{Shifted Schur measures}
%

We recall the Schur $Q$-function and the shifted Schur measure.
The following facts are known in \cite[III-8]{Macdonald} and \cite{TracyWidom200?}. 

A non-increasing sequence  $\lambda= (\lambda_1, \lambda_2, \dots)$ of non-negative integers is called a partition of $N$ 
if the size $|\lambda|:= \sum_{j \geq 1} \lambda_j$ equals $N$.
We denote the number of non-zero parts of $\lambda$ by $\ell(\lambda)$ and we call it the length of $\lambda$.
A partition $\lambda$ is called {\it strict} if and only if all parts of $\lambda$ are distinct
and then we write $\lambda \vDash N$.
Let $\mD_N$ be the set of all strict partitions of $N$ and $\mD$ the set of all strict partitions, 
i.e., $\mD= \cup_{N=0}^{\infty} \mD_N$.  

Let $X=(X_1,X_2, \dots)$ and $Y=(Y_1,Y_2,\dots)$ be infinite many variables.
The symmetric functions $q_n(X)$ $(n \geq 0)$ are defined via the generating function
  \begin{equation*}
   Q(z)=Q_X(z) = \prod_{i=1}^{\infty} \frac{1+X_i z}{1- X_i z} = \sum_{n=0}^{\infty} q_n (X) z^n.
  \end{equation*}
In particular, we have $q_0=1$.
Since
$$
  \log \prod_{i=1}^{\infty} \frac{1+X_i z}{1- X_i z} 
  = \sum_{i=1}^\infty \sum_{n=1}^{\infty} \frac{1- (-1)^n}{n}X_i^n z^n
  = \sum_{n=1,3,5,\dots} \frac{2}{n} p_n(X) z^n,
$$ 
where $p_n(X)=\sum_{i=1}^\infty X_i^n$ is the power-sum function,
the function $Q(z)$ is also expressed as 
\begin{equation} \label{Q(z)2}
Q(z)= \exp \(\sum_{n=1,3,5,\dots} \frac{2}{n} p_n(X) z^n \).
\end{equation}

For $\lambda=(\lambda_1,\lambda_2,\dots) \in \mD$ of length $\leq m$,
the Schur $Q$-function $Q_{\lambda}(X)$ is defined as 
the coefficient of $z^{\lambda}=z_1^{\lambda_1} z_2^{\lambda_2} \cdots z_m^{\lambda_m}$ in
  \begin{equation} \label{generatingfunction}
   Q(z_1, z_2, \dots, z_m) = \prod_{i=1}^m Q(z_i) \prod_{1 \leq i<j \leq m} \frac{z_i-z_j}{z_i+z_j}.
  \end{equation} 

For $r>s \geq 0$, we define
$$
Q_{(r,s)} = q_r q_s + 2 \sum_{i=1}^s (-1)^i q_{r+i} q_{s-i}
$$
and $Q_{(r,s)} = - Q_{(s,r)}$ for $r \leq s$.
We may write $\lambda$ in the form
$\lambda = (\lambda_1, \lambda_2, \dots, \lambda_{2n})$
where $\lambda_1 > \lambda_2 > \dots > \lambda_{2n} \geq 0$.
Then the $2n \times 2n$ matrix 
$$
M_{\lambda} = ( Q_{(\lambda_i, \lambda_j)})_{1 \leq i,j \leq 2n}
$$
is skew symmetric, and the Schur $Q$-function $Q_{\lambda}$ is also given by
\begin{equation}
  Q_{\lambda} = \Pf (M_{\lambda}),
\end{equation}
where $\Pf$ stands for the pfaffian.
The Schur $P$-function $P_{\lambda}$ is defined by
$P_{\lambda} = 2^{-\ell(\lambda)} Q_{\lambda}$.

The {\it shifted Schur measure} is a (formal) probability measure on $\mD$ defined by
  \begin{equation} \label{ShiftedSchur}
   \PP_{\mathrm{SS}} ( \{ \lambda \} ) = \frac{1}{ Z_{\mathrm{SS}} } Q_{\lambda} (X) P_{\lambda} (Y)
  \end{equation}
for each $\lambda \in \mD$. 
Here the normalization constant $Z_{\mathrm{SS}}$ is determined by
  $$
   Z_{\mathrm{SS}} = \sum_{\lambda \in \mD} Q_{\lambda}(X) P_{\lambda}(Y) = 
    \prod_{i,j=1}^{\infty} \frac{1 + X_i Y_j}{1 - X_i Y_j},
  $$
where the second equality is the Cauchy identity for Schur $Q$-functions (\cite[p.255]{Macdonald}).
Further, from \eqref{Q(z)2},
the constant $Z_{\mathrm{SS}}$ is also expressed as
\begin{equation} \label{Z_SS}
Z_{\mathrm{SS}} = \exp \( \sum_{n=1,3,5,\dots} \frac{2}{n} p_n(X) p_n(Y) \).
\end{equation}

%
\section{Correlation functions of the shifted Schur measure}
%

In this section, we prove the main theorem.
We identify each strict partition $\lambda=(\lambda_1, \lambda_2, \dots, \lambda_{\ell})$
$(\lambda_1 > \lambda_2 > \dots > \lambda_{\ell} >0)$ 
with the finite set $\{ \lambda_1, \lambda_2, \dots, \lambda_{\ell} \}$
of positive integers.
Define
the correlation function of the shifted Schur measure $\PP_{\mathrm{SS}}$ by 
  \begin{equation}
   \rho_{\mathrm{SS}}(A) := \PP_{\mathrm{SS}} ( \{ \lambda \in \mD | \lambda \supset A \}) = 
   \frac{1}{Z_{\mathrm{SS}}} \sum_{\lambda \supset A} Q_{\lambda}(X) P_{\lambda} (Y)
  \end{equation}
for a finite subset $A \subset \bZ_{>0}$.
The function $\rho_{\mathrm{SS}}(A)$ has a pfaffian expression.

\noindent
\begin{thm} \label{corSS}
For a finite subset $A= \{k_1, \dots, k_N \} \subset \bZ_{>0}$, we have
   \begin{equation}
    \rho_{\mathrm{SS}} (A) =\Pf (M(A)_{i,j})_{1 \leq i < j \leq 2N },
   \end{equation}
where the entry $M(A)_{i,j}$ of the skew symmetric matrix $M(A)$ is given by
  $$
   M(A)_{i,j} = \begin{cases}
                \bK(k_i, k_j ), & \text{for} \quad 1 \leq i<j \leq N, \\
                \bK(k_i, -k_{2N-j+1}), & \text{for} \quad 1 \leq i \leq N < j \leq 2N, \\
                \bK(-k_{2N-i+1}, -k_{2N-j+1}), & \text{for} \quad N< i<j \leq 2N,
                \end{cases}
  $$
and $\bK (u,v)$ is defined as $\epsilon (u,v)$ times the coefficient of $z^u w^v$ in the formal series
    $$
      \frac{1}{2} \bJ(z;X,Y) \bJ(w;X,Y) \frac{z-w}{z+w}.
    $$
Here 
$\bJ (z;X,Y)$ is defined by
  \begin{equation} \label{bJ}
    \bJ (z;X,Y) := Q_X(z) Q_Y (-z^{-1}) 
     = \prod_{i=1}^{\infty} \frac{1+X_i z}{1-X_i z} \frac{1 -Y_i z^{-1}}{1+ Y_i z^{-1}}
  \end{equation}
and $\epsilon (u,v)$ is given by
  \begin{equation} \label{epsilon}
    \epsilon(u,v) = \begin{cases}
                     1, & \text{for} \quad u,v >0, \\                    
                     (-1)^v, & \text{for} \quad u>0, v<0, \\
                     (-1)^{u+v}, & \text{for} \quad u,v <0.
                    \end{cases}
  \end{equation}
\end{thm}

\medskip

\begin{remark}
The correlation function of the Schur measure
is given by a determinant, see Theorem 1 in \cite{Okounkov2001}.
\qed
\end{remark}

\medskip

We prove Theorem \ref{corSS} by employing the exterior algebra.
Let $V$ be a module on $\bZ[X_1,X_2, \dots, Y_1,Y_2, \dots]$ spanned by $\be_k \ (k =1,2,\dots)$.
The exterior algebra $\bigwedge V$ is spanned by vectors
  $$
    \bv_{\lambda}= \be_{\lambda_1} \wedge \be_{\lambda_2} \wedge \dots \wedge \be_{\lambda_{\ell}},
  $$
where $\lambda=(\lambda_1, \dots, \lambda_{\ell}) \in \mD$ ($\lambda_1 > \dots > \lambda_{\ell} \geq 1$).
In particular, we have $\bv_{\emptyset} = 1$.
We give $\bigwedge V$ the inner product
  $$
    \langle \bv_{\lambda}, \bv_{\mu} \rangle = \delta_{\lambda, \mu} 2^{-\ell(\lambda)}.
  $$
Putting $\be^{\vee}_k = 2 \be_k$ and $\bv^{\vee}_{\lambda}= 
\be^{\vee}_{\lambda_1} \wedge \dots \wedge \be^{\vee}_{\lambda_{\ell}}= 2^{\ell} \bv_{\lambda}$,
the bases $(\bv_{\lambda})_{\lambda \in \mD}$ and $(\bv^{\vee}_{\lambda})_{\lambda \in \mD}$
are dual to each other.

We define the operator $\psi_k$ ($k \geq 1$) on $\bigwedge V$ by
  $$
    \psi_k \bv_{\lambda} = \be_k \wedge \bv_{\lambda}
  $$
and let $\psi^*_k$ be the adjoint operator of $\psi_k$ with respect to the inner product defined above.
The operator $\psi^*_k$ is then explicitly given by
  $$
    \psi^*_k \bv_{\lambda} = \sum_{i=1}^{\ell(\lambda)} \frac{(-1)^{i-1}}{2} \delta_{k, \lambda_i} 
    \be_{\lambda_1} \wedge \dots \wedge \widehat{\be_{\lambda_i}} \wedge \dots \wedge \be_{\lambda_{\ell}}.
  $$
These operators satisfy the following commutation relations
  \begin{equation}
    \begin{array}{lcr} \label{cr1}
    \psi_i \psi^*_j + \psi^*_j \psi_i = \delta_{i,j} \frac{1}{2}, &
    \psi_i \psi_j = - \psi_j \psi_i, &
    \psi^*_i \psi^*_j = -\psi^*_j \psi^*_i. 
   \end{array}
  \end{equation}

Since
  \begin{equation} \label{corr}
    \psi_k \psi^*_k \bv_{\lambda} = \begin{cases}
                            \frac{1}{2} \bv_{\lambda}, & \text{if} \ k \in \lambda, \\
                            0,                         & \text{otherwise},
                           \end{cases} 
  \end{equation}
we see that $\(\prod_{k \in A} 2\psi_k \psi^*_k\) \bv_{\lambda}$ is equal 
to $\bv_{\lambda}$ if $A \subset \lambda$ and to $0$ otherwise.

Define the self-adjoint operator $S$ by $S \bv_{\lambda}= (-1)^{\ell(\lambda)} \bv_{\lambda}$ for any $\lambda \in \mD$.
The operators satisfy the relations 
 \begin{equation} \label{cr2}
  \begin{array}{lcr}
    S^2=1,        &  \psi_k S= -S \psi_k,   &    \psi^*_k S = - S \psi^*_k.
  \end{array}
 \end{equation}

For each odd positive integer $n$, we define the operators $\alpha_n$ and $\alpha_{-n}$ by
  \begin{align*}
   \alpha_n &:= 2 \sum_{j=1}^{\infty} \psi_j \psi^*_{n+j} 
     + S \psi^*_n + 2 \sum_{j=1}^{\frac{n-1}{2}} (-1)^j \psi^*_j \psi^*_{n-j}, \\
   \alpha_{-n} := \alpha^*_n &= 2 \sum_{j=1}^{\infty} \psi_{n+j} \psi^*_j + \psi_n S + 2 \sum_{j=1}^{\frac{n-1}{2}}
                  (-1)^j \psi_{n-j} \psi_j.
 \end{align*} 
It follows from \eqref{cr1} and \eqref{cr2} that
  \begin{equation} \label{bracketalpha}
    [ \alpha_n, \alpha_m] = \frac{n}{2} \delta_{n,-m}
  \end{equation}
for any odd integers $n$ and $m$, where $[ \ , \ ]$ is the commutator; $[a,b]=ab-ba$.

If we put
\begin{equation}
   \tpsi_k= \begin{cases} 
             \psi_k, & \text{for} \quad k \geq 1, \\
             S/2,     & \text{for} \quad k=0, \\
             (-1)^k \psi^*_{-k}, & \text{for} \quad k \leq -1,
            \end{cases}
  \end{equation}
and $\psi(z)= \sum_{k \in \bZ} z^k \tpsi_k$,
then by \eqref{cr1} and \eqref{cr2} we see that
  \begin{equation} \label{bracketpsi}
      [ \alpha_n, \psi(z) ]= z^n \psi(z) \quad \text{for any odd integer $n$}
  \end{equation}
and 
  \begin{equation} \label{cr3}
   \langle \tpsi_k \tpsi_l \bv_{\emptyset}, \bv_{\emptyset} \rangle = 0 \qquad \text{unless} \qquad l=-k \geq 0.
  \end{equation}
It follows from \eqref{cr3} that
 \begin{equation} \label{psi(z)psi(w)}
   \begin{split}
    \langle \psi(z) \psi(w) \bv_{\emptyset}, \bv_{\emptyset} \rangle
   &= \left\langle \( \frac{S^2}{4}+ \sum_{k \geq 1} (-1)^k z^{-k} w^k \psi^*_k \psi_k \) \bv_{\emptyset}, 
   \bv_{\emptyset} \right\rangle \\
   &=\frac{1}{4} + \sum_{k \geq 1} \frac{1}{2} \( -\frac{w}{z}\)^k = \frac{z-w}{4(z+w)}.
    \end{split}
  \end{equation}
Note that the operator $\alpha_n$ is expressed as $\alpha_n = \sum_{k \in \bZ} (-1)^k\tpsi_{k-n} \tpsi_{-k}$.

Put 
 \begin{equation*}
   \Gamma_{\pm} (X) = \exp \( \sum_{n=1,3,5,\dots} \frac{2 p_n(X)}{n} \alpha_{\pm n} \).
 \end{equation*} 
Observe that
  \begin{align}
   \Gamma_{+} \bv_{\emptyset} &= \bv_{\emptyset}, \label{Gamma+} \\
   \Gamma^*_{\pm} &= \Gamma_{\mp}, \label{Gammapm} \\
   \Gamma_{+}(X) \Gamma_{-}(Y) 
    &= Z_{\mathrm{SS}} \Gamma_{-}(Y) \Gamma_{+} (X). \label{Gammacom}
  \end{align}
The equality \eqref{Gammacom} is obtained from \eqref{bracketalpha} and \eqref{Z_SS}.
By \eqref{Q(z)2} and \eqref{bracketpsi}, we have
 \begin{equation} \label{Gammapsi}
   \Gamma_{\pm}(X) \psi(z) = Q_X(z^{\pm 1}) \psi(z) \Gamma_{\pm}(X).  
 \end{equation}

The Schur $Q$-function is given as a matrix element of $\Gamma_{-}$ as follows.

\begin{prop} \label{propQfunction}
For each $\lambda \in \mD$, we have
  \begin{equation} \label{innerproduct}
    \langle \Gamma_{-} (X) \bv_{\emptyset}, \bv^{\vee}_{\lambda} \rangle = Q_{\lambda}(X).
  \end{equation}
More generally, for $\lambda,\mu \in \mD$,
\begin{equation} \label{innerproduct2}
\langle \Gamma_{-} (X) \bv_{\mu}, \bv^{\vee}_{\lambda} \rangle = Q_{\lambda/\mu}(X),
\end{equation}
where $Q_{\lambda/\mu}(X)$ is a skew Hall-Littlewood function.
\end{prop}

\begin{proof}
Write $\lambda$ in the form  $\lambda_1 > \lambda_2 > \cdots > \lambda_{2n} \geq 0$.
Since $\bv^{\vee}_{\lambda} = 2^{2n} \tpsi_{\lambda_1} \cdots \tpsi_{\lambda_{2n}} \bv_{\emptyset}$, 
the left hand side in \eqref{innerproduct} is equal to the coefficient of $z_1^{\lambda_1} \cdots z_{2n}^{\lambda_{2n}}$ in
the expansion of
  \begin{equation} \label{ip}
    2^{2n} \langle \Gamma_{-} (X) \bv_{\emptyset}, \psi(z_1) \cdots \psi(z_{2n}) \bv_{\emptyset} \rangle.
  \end{equation}
It follows from \eqref{Gamma+}, \eqref{Gammapm} and \eqref{Gammapsi} that \eqref{ip} equals 
  $$
   2^{2n} Q(z_1) \cdots Q(z_{2n}) \langle \psi(z_1) \cdots \psi(z_{2n}) \bv_{\emptyset}, \bv_{\emptyset} \rangle.
  $$
By \eqref{generatingfunction},
in order to prove \eqref{innerproduct} it is sufficient to show 
  $$
    2^{2n} \langle \psi(z_1) \cdots \psi(z_{2n}) \bv_{\emptyset}, \bv_{\emptyset} \rangle 
     = \Pf \( \frac{z_i-z_j}{z_i+z_j} \) = \prod_{1 \leq i<j \leq 2n} \frac{z_i-z_j}{z_i+z_j}.
  $$
Note the second equality is well-known (see e.g. \cite[III-8, Ex.5]{Macdonald}). 
From \eqref{cr1}, \eqref{cr2} and \eqref{cr3},
we see that
\begin{align*}
   & 2^{2n} \langle \psi(z_1) \cdots \psi(z_{2n}) \bv_{\emptyset}, \bv_{\emptyset} \rangle \\
     =& \sum_{k=2}^{2n} (-1)^k 4 \langle \psi(z_1) \psi(z_k) \bv_{\emptyset}, \bv_{\emptyset} \rangle         
     4^{n-1} \langle \psi(z_2) \cdots \widehat{\psi(z_k)} \cdots \psi(z_{2n}) \bv_{\emptyset}, \bv_{\emptyset} \rangle.
\end{align*}
Therefore, by the expansion formula of a pfaffian, we obtain
  $$
   2^{2n} \langle \psi(z_1) \cdots \psi(z_{2n}) \bv_{\emptyset}, \bv_{\emptyset} \rangle =
    \Pf ( 4 \langle \psi(z_i) \psi(z_j) \bv_{\emptyset}, \bv_{\emptyset} \rangle ).
  $$
Hence the claim follows from \eqref{psi(z)psi(w)}.
The generating function of $Q_{\lambda/\mu}$ in \cite[III-8, Ex.9]{Macdonald}
yields the second formula \eqref{innerproduct2} by a discussion similar to the above.
\end{proof}

From \eqref{corr}, \eqref{Gammapm} and \eqref{innerproduct}, 
the correlation function is expressed as
  \begin{equation*}
    \rho_{\mathrm{SS}} (A) = \frac{1}{Z_{\mathrm{SS}}} \sum_{ \lambda \supset A} 
                                                      Q_{\lambda} (X) P_{\lambda}(Y)
   =  \frac{1}{Z_{\mathrm{SS}}} \left\langle \Gamma_{+} (X) \( \prod_{k \in A} 2 \psi_k \psi^*_k \) \Gamma_{-} (Y)
       \bv_{\emptyset}, \bv_{\emptyset} \right\rangle.
  \end{equation*}
It follows from \eqref{Gamma+}, \eqref{Gammapm} and \eqref{Gammacom} that
  \begin{equation} \label{entry}
   \rho_{\mathrm{SS}} (A) = \left\langle \( \prod_{k \in A} 2 \Psi_k \Psi^*_k \) 
                               \bv_{\emptyset}, \bv_{\emptyset} \right\rangle,
  \end{equation}
where we put
  \begin{equation} \label{Psi}
   \begin{array}{ccc}
    \Psi_k = \Ad (G) \psi_k, & \Psi^*_k = \Ad (G) \psi^*_k, & G= \Gamma_{+}(X) \Gamma_{-}(Y)^{-1}.
   \end{array}
  \end{equation}
Using \eqref{Gammapsi}, we have
  \begin{equation} \label{Ad(G)psi(z)}
    \Ad(G) \psi(z)  = \bJ (z;X,Y) \psi(z),
  \end{equation} 
where $\bJ (z;X,Y)$ is defined in \eqref{bJ}.

\begin{lem} \label{LemCor}
 We have
  \begin{equation*} 
   \rho_{\mathrm{SS}} (A) = 
   \Pf ( \widetilde{M}(A)_{i,j})_{1 \leq i<j \leq 2N}.
  \end{equation*}
Here the entry of the skew symmetric matrix $\widetilde{M}(A)$ is given by
  \begin{equation} \label{pf}
    \widetilde{M}(A)_{i,j} = \begin{cases}
      2\langle \Psi_{k_i} \Psi_{k_j} \bv_{\emptyset}, \bv_{\emptyset} \rangle, & \text{for} \quad 1 \leq i<j \leq N, \\
      2\langle \Psi_{k_i} \Psi^{*}_{k_{2N-j+1}} \bv_{\emptyset}, \bv_{\emptyset} \rangle,  
                                                                    & \text{for} \quad 1 \leq i \leq N < j \leq 2N, \\
             2\langle \Psi^{*}_{k_{2N-i+1}} \Psi^{*}_{k_{2N-j+1}} \bv_{\emptyset}, \bv_{\emptyset} \rangle,   
                                                                            & \text{for} \quad N< i<j \leq 2N.
                         \end{cases}
  \end{equation}
\end{lem}

\begin{proof} 
From \eqref{cr1} and \eqref{Psi}, we have $\Psi_k \Psi^*_l = - \Psi^*_l \Psi_k$ ($k \not= l$).
Therefore we obtain
  $$
    \left\langle \( \prod_{k \in A}  \Psi_k \Psi^*_k \) 
                               \bv_{\emptyset}, \bv_{\emptyset} \right\rangle =
    \langle \Psi_{k_1} \Psi_{k_2} \cdots \Psi_{k_N} \Psi^*_{k_N} \cdots \Psi^*_{k_1} \bv_{\emptyset}, 
     \bv_{\emptyset} \rangle.
  $$
By \eqref{Ad(G)psi(z)},
the operator $\Psi_j$ and $\Psi^*_j$, respectively, is expressed as a linear combination of $\tpsi_n$'s over 
$\bZ[X_1,X_2, \dots, Y_1,Y_2, \dots ]$.
Hence if we abbreviate 
$\tPsi_j = \Psi_{k_j}$ for $1 \leq j \leq N$ and $\tPsi_j= \Psi^*_{k_{2N-j+1}}$ for $N+1 \leq j \leq 2N$,
we have
  $$
  \langle \Psi_1 \Psi_2 \cdots \Psi_N \Psi^*_N \cdots \Psi^*_1 \bv_{\emptyset}, \bv_{\emptyset} \rangle
   = \Pf ( \langle \tPsi_i \tPsi_j \bv_{\emptyset}, \bv_{\emptyset} \rangle)_{1 \leq i< j \leq 2N}
  $$ 
by a discussion similar to the proof of Proposition \ref{propQfunction}.
Thus, by \eqref{entry}, we obtain the lemma.
\end{proof}

\begin{proof}[Proof of Theorem \ref{corSS}]
We compute entries in the right hand side of \eqref{pf}.
It follows from \eqref{psi(z)psi(w)} and \eqref{Ad(G)psi(z)} that
  $$
   \langle 2 \Psi (z) \Psi(w) \bv_{\emptyset}, \bv_{\emptyset} \rangle 
  =  \frac{1}{2} \bJ(z;X,Y) \bJ(w;X,Y) \frac{z-w}{z+w},
  $$
where $\Psi(z) = \Ad (G) \psi(z)$.
Since the coefficient of $z^k$ ($k \in \bZ \setminus \{0 \}$) in $\Psi(z)$ is equal to $\Psi_k$ if $k >0$ and
to $(-1)^k \Psi^*_{-k}$ if $k<0$,
we can easily see the theorem from Lemma \ref{LemCor}.
\end{proof}

\begin{remark}
Though Jing \cite{Jing} obtains the expression of Schur $Q$-functions by vertex operators 
with the commutator relation \eqref{bracketalpha}
it seems very hard to obtain the result in Theorem \ref{corSS} using these vertex operators.
\qed
\end{remark}

%
\section{Applications} \label{Applications}
%

As an application of Theorem \ref{corSS}, we give a limit distribution of $\lambda_j$'s
with respect to a specialization of the shifted Schur measure.

%
\subsection{A shifted version of the Plancherel measure} 
%

We define a measure similar to the Plancherel measure on $\mP_N$
by means of the shifted Robinson-Schensted-Knuth (RSK) correspondence (see e.g. \cite{HH}).

A {\it shifted shape} $\Sh (\lambda)$ associated with a strict partition $\lambda$
is obtained by replacing the $i$-th row to the right by $i-1$ boxes for $i \geq 1$ from the Young diagram $\lambda$. 
A {\it standard shifted tableau} $T$ of the shifted shape $\lambda \vDash N$ 
is an assignment of $1,2, \dots, N$ to each box in 
the shifted shape $\Sh (\lambda)$ such that 
entries in $T$ are increasing across rows and down columns.
For example, 
$$
\begin{matrix} 1 & 2 & 4 & 6 \\
                 & 3 & 5 & 8 \\
                 &   & 7 &   \end{matrix}
$$
is a standard shifted tableau of shape $\lambda=(4,3,1)$.

Let $g^{\lambda}$ be the number of standard shifted tableaux of shape $\lambda$.
It is known that 
$g^{\lambda}$ is explicitly given by
$$
g^{\lambda}= \frac{|\lambda|!}{\lambda_1 ! \lambda_2 ! \cdots \lambda_{\ell} !} 
\prod_{1 \leq i<j \leq \ell} \frac{\lambda_i-\lambda_j}{\lambda_i+\lambda_j}
$$
(see e.g. \cite[III-8, Ex.12]{Macdonald}).
By means of the shifted RSK we can see that
\begin{equation} \label{eq.shifted}
\sum_{\lambda \vDash N} 2^{N- \ell(\lambda)} (g^{\lambda})^2 = N!
\end{equation} 
(see \cite{HH}).

In view of the equality \eqref{eq.shifted},
we define a probability measure on $\mD_N$, that is, 
we assign to each $\lambda \in \mD_N$ the probability
\begin{equation} \label{SPl}
   \PP_{\mathrm{SPl},N} (\{ \lambda \}) = \frac{2^{N- \ell (\lambda)} }{N!} (g^{\lambda})^2.
\end{equation}
This measure, which is noted in \cite{TracyWidom200?},
can be regarded as a shifted version of the Plancherel measure defined in \eqref{Plancherel} in a combinatorial sense.

%
\subsection{Ascent pairs for a permutation}
%

The measure defined in \eqref{SPl} is related to the so-called ascent pair for a permutation.
For $\pi = ( \pi(1), \pi(2), \dots, \pi(N)) \in \mS_N$,
an {\it ascent pair} for $\pi$ is a pair  $(\phi^{\mathrm{de}}, \phi^{\mathrm{in}})$ 
of a decreasing subsequence $\phi^{\mathrm{de}}=(\pi(i_1)> \dots> \pi(i_k)), \ i_1< \dots <i_k$ and 
an increasing subsequence $\phi^{\mathrm{in}}=(\pi(j_1)< \dots< \pi(j_l)), \ j_1< \dots <j_l$ 
of $\pi$ such that the sequence
$$
(\pi(i_k), \dots, \pi(i_1), \pi(j_1), \dots, \pi(j_l))
$$
is weakly increasing (i.e. the inequality $\pi(i_1) \leq \pi(j_1)$ is satisfied).
We define the length of the ascent pair  $(\phi^{\mathrm{de}}, \phi^{\mathrm{in}})$ by $k+l-1$.
Denote the length of the longest ascent pair for $\pi$ by $L(\pi)$.

\begin{example}
For a permutation
$$
\pi=\begin{pmatrix} 1 & 2 & 3 & 4 & 5 & 6 & 7 & 8 & 9 \\ 4 & 7 & 1 & 9 & 6 & 3 & 5 & 8 & 2 \end{pmatrix}
$$
the pair $( \phi^{\mathrm{de}}, \phi^{\mathrm{in}})$, where
$\phi^{\mathrm{de}}= (4,3,2)$ and $\phi^{\mathrm{in}} = (4,7,9)$,
is the ascent pair with length 5.
Since this is the longest ascent pair for $\pi$,
we have $L(\pi)=5$.
\qed
\end{example}

\medskip

By the shifted RSK,
the distribution of $L(\pi)$ with respect to the uniform measure on $\mS_N$
equals the distribution of $\lambda_1$ with respect to the measure $\PP_{\mathrm{SPl},N}$ on $\mD_N$, i.e.,
\begin{equation} \label{PlancherelAscent}
\PP_{\mathrm{uniform},N} (\{ \pi \in \mS_N | L(\pi) = h \}) = \PP_{\mathrm{SPl},N} 
( \{ \lambda \in \mD_N | \lambda_1 = h \}).
\end{equation}

\subsection{Limit distributions}

We consider the random point process on $\bR$ (see the Appendix in \cite{BOO}) whose correlation functions
$\rho_{\Airy}(X) = \PP_{\Airy} ( \{ Y \subset \bR \ | \ \# Y < \infty, \ X \subset Y \} )$ 
for any finite subset $X=\{x_1, \dots, x_k \} \subset \bR$
are given by $\rho_{\Airy}(X) = \det (K_{\Airy}(x_i,x_j))_{1 \leq i,j \leq k}$.
Here $K_{\Airy}$ is the Airy kernel defined in \eqref{Airykernel}.
Let 
$\zeta = (\zeta_1 > \zeta_2 > \cdots ) \in \bR^{\infty}$
be its random configuration.
The random variables $\zeta_i$'s are called the {\it Airy ensemble}.
It is known that the Airy ensemble describes the behavior of 
the largest eigenvalue of a GUE matrix, the 2nd largest one, and so on, see \cite{TracyWidom1994}.

Theorem 4 in \cite{BOO} (see also \cite{JohanssonDis,Okounkov2000}) asserts that  
the random variables 
\begin{equation} \label{randomvariables}
\frac{\lambda_i - 2 \sqrt{N}}{N^{1/6}}, \quad i=1,2,\dots, \quad \lambda=(\lambda_1,\lambda_2, \dots) 
\in \mP_N \end{equation}
with respect to the Plancherel measure defined by \eqref{Plancherel}
converge, in the joint distribution, 
to the Airy ensemble as $N \to \infty$.
The following theorem is a shifted version of this result.

\begin{thm} \label{mainthm}
As $N \to \infty$, the random variables 
\begin{equation} \label{randomvariables2}
\frac{\lambda_i - 2 \sqrt{2N}}{(2N)^{1/6}}, \quad i=1,2,\dots
\end{equation}
with respect to the measure $\PP_{\mathrm{SPl},N}$ on $\mD_N$
converge to the Airy ensemble, in joint distributions.
\end{thm}

Compare \eqref{randomvariables2} with \eqref{randomvariables}.
Especially, since the distribution of $\zeta_1$ in the Airy ensemble is given by the Tracy-Widom distribution,
we immediately see the following result from \eqref{PlancherelAscent}.

\begin{cor} \label{maincor}
We have
$$
\lim_{N \to \infty} \PP_{\mathrm{uniform},N} \( \frac{L - 2 \sqrt{2N}}{(2N)^{1/6}} <s \) = F_2(s).
$$ 
\end{cor}

Compare with \eqref{thmBDJ}.
Theorem \ref{mainthm} is proved by computing the correlation function 
of the so-called {\it poissonization} of the measure $\PP_{\mathrm{SPl},N}$.
For $\xi >0$, 
we define the poissonization $\PP^{\xi}_{\mathrm{PSP}}$ of the measure $\PP_{\mathrm{SPl},N}$ by
\begin{equation} \label{PSP}
  \PP^{\xi}_{\mathrm{PSP}} (\{ \lambda \}) = e^{-\xi} \sum_{N=0}^{\infty} \frac{\xi^N}{N!} 
  \PP_{\mathrm{SPl},N} (\{ \lambda \} )
    = e^{-\xi} \xi^{|\lambda|} 2^{|\lambda| - \ell (\lambda) } \( \frac{ g^{\lambda} }{ | \lambda | !} \)^2 
\end{equation}
for $\lambda \in \mD$.
Here $\PP_{\mathrm{SPl},N}(\{ \lambda \}) = 0$ unless $\lambda \vDash N$.
Then we have the

\begin{thm} \label{thmpoi}
For any fixed $M \geq 1$ and any $a_1, \dots, a_m \in \bR$ we have
\begin{equation}
\lim_{\xi \to \infty} \PP^{\xi}_{\mathrm{PSP}} 
\( \left\{ \lambda \in \mD \Biggm| \frac{\lambda_i - 2 \sqrt{2 \xi} }{(2 \xi)^{\frac{1}{6}}}
< a_i,\  1 \leq i \leq M \right\} \) 
= \PP_{\Airy} ( \zeta_i < a_i,\  1 \leq i \leq M),
\end{equation}
where $\zeta_1 > \zeta_2 > \cdots$ is the Airy ensemble.
\end{thm}

Since Theorem \ref{mainthm} can be proved from Theorem \ref{thmpoi} 
by using the {\it depoissonization technique} developed in \cite{Johansson1998},
we omit the proof, see \cite{BOO}.

%
\subsection{The proof of Theorem \ref{thmpoi}}
%

The measure $\PP_{\mathrm{PSP}}^\xi$ can be obtained by a specialization of the shifted Schur measure.
Actually,
since the Schur $Q$-function can be expanded as (see \cite{Macdonald})
  $$
   Q_{\lambda}(X) = \sum_{\rho = 1^{m_1} 3^{m_3} \cdots} 2^{ \ell (\rho) } X_{\rho}^{\lambda} (-1) 
                      \prod_{i : \odd} \frac{ p_k (X)^{m_i} }{ m_i ! i^{m_i}},
  $$
where $X_{(1^{|\lambda|})}^{\lambda} (-1)= g^{\lambda}$ 
(see \cite[III-8, Ex.12]{Macdonald}),
if we make a specialization such as $p_k(X)= p_k(Y) = \sqrt{ \frac{\xi}{2} } \delta_{k1} \ (k \geq 1)$, 
then we have
  $$
   Q_{\lambda} = (2 \xi)^{ \frac{| \lambda |}{2} }  \frac{g^{\lambda}}{| \lambda | !}.
  $$ 
Hence the shifted Schur measure in \eqref{ShiftedSchur} becomes 
the measure $\PP^{\xi}_{\mathrm{PSP}}$ in \eqref{PSP}.
 
Let $\rho_{\mathrm{PSP}}^\xi$ be the correlation function of the measure $\PP_{\mathrm{PSP}}^\xi$.

\begin{prop} \label{propAiry}
We have
$$
 \lim_{\xi \to +\infty} (2\xi)^{N/6} \rho^{\xi}_{\mathrm{PSP}} 
( \{ [2\sqrt{2\xi} +(2\xi)^{1/6} x_1], \dots, [2\sqrt{2\xi} +(2\xi)^{1/6} x_N] \})
= \det (K_{\Airy} (x_i,x_j))_{1 \leq i,j \leq N}.
$$
The limit is uniform for $(x_1, \dots, x_N)$ on a compact set of $\bR^N$.
\end{prop}

This proposition follows immediately
from Theorem \ref{corSS} and the following lemma.

\begin{lem} \label{lemlimit}
We have
\begin{align}
 (2\xi)^{\frac{1}{6}} \bK_{\mathrm{B}} 
(2 \sqrt{2\xi} + x (2\xi)^{\frac{1}{6}}, 2 \sqrt{2\xi} + y (2\xi)^{\frac{1}{6}} )  &\to 0 \label{bK++}, \\ 
 (2\xi)^{\frac{1}{6}} \bK_{\mathrm{B}} 
(2 \sqrt{2\xi} + x (2\xi)^{\frac{1}{6}}, -(2 \sqrt{2\xi} + y (2\xi)^{\frac{1}{6}}) )  &\to K_{\Airy} (x,y), \label{bK+-} \\
 (2\xi)^{\frac{1}{6}} \bK_{\mathrm{B}} 
(-(2 \sqrt{2\xi} + x (2\xi)^{\frac{1}{6}}), -(2 \sqrt{2\xi} + y (2\xi)^{\frac{1}{6}}) )  &\to 0, \label{bK--}
\end{align}
as $\xi \to \infty$, uniformly in $x$ and $y$ on compact sets in $\bR$.
\end{lem}

\begin{proof}
By the specialization $p_k(X)=p_k(Y) = \sqrt{\xi/2} \delta_{k1}$, 
the function $\bJ(z;X,Y)$ in \eqref{bJ} becomes $e^{\sqrt{2\xi} (z-z^{-1})}$,
which is the generating function of Bessel functions.
Therefore, in order to prove Lemma \ref{lemlimit}, we evaluate integrals of the form
$$
\( \frac1{2 \pi \sqrt{-1}} \) \iint e^{2\xi (z-z^{-1} +w -w^{-1})} \frac{z-w}{z+w} \frac{\d z \d w}{z^{u+1} w^{v+1}},
$$
where the contours are two unit circles and  $u = \pm (2 \sqrt{2\xi} + x (2\xi)^{ \frac{1}{6} } )$ and 
$v = \pm (2 \sqrt{2\xi} + y (2\xi)^{ \frac{1}{6} } )$.
Then Lemma \ref{lemlimit} is obtained by a similar discussion in \cite{TracyWidom200?}.
We leave the detail for readers.
\end{proof}

\begin{proof}[Proof of Theorem \ref{thmpoi}]
The proof follows from Proposition \ref{propAiry} and the discussion in \cite{BOO}.  
\end{proof}

%
\subsection{The $\alpha$-specialized shifted Schur measure}
%

Let $\alpha$ be a real number such that $0 < \alpha <1$ and
let $m$ and $n$ be positive integers.
We put $X_i = Y_j = \alpha$ for $1 \leq i \leq m$ and $1 \leq j \leq n$,
and let the rest be zero
in the definition of the shifted Schur measure.
This is called the {\it $\alpha$-specialization},
see \cite{TracyWidom200?} and \cite{Matsumoto}.
Using Theorem \ref{corSS},
we also give a limit distribution of $\lambda_i$'s with respect to the $\alpha$-specialized shifted Schur measure.
Denote by $\PP_{\mathrm{SS},\sigma}$ the $\alpha$-specialized shifted Schur measure, 
where $\sigma= (m,n, \alpha)$ denotes the set of parameters above, and
put $\tau=m/n$.

\begin{thm} \label{alphashifted}
There exist positive constants $c_1= c_1 (\alpha,\tau)$ and $c_2 =c_2(\alpha,\tau)$ such that
$$
\lim_{n \to \infty} \PP_{\mathrm{SS},\sigma} \( \left\{ \lambda \in \mD \  \big| \ 
\frac{\lambda_i - c_1 n}{c_2 n^{1/3}} < a_i, \ 1 \leq i \leq M \right\} \) 
= \PP_{\Airy} ( \zeta_i < a_i, \  1 \leq i \leq M)
$$
holds for any $M \geq 1$ and any $a_1, \dots, a_M \in \bR$.
\end{thm}

When $M=1$,
this theorem gives the result in \cite{TracyWidom200?}.
Although they assume that $\alpha$ and $\tau$ satisfy the relation $\alpha^2 < \tau < \alpha^{-2}$, 
we can remove this assumption as they expect in the footnote 3 of that paper.

We write $\rho_{\mathrm{SS}}$, $M$, $\bJ(z;X,Y)$ and $\bK$ in Theorem \ref{corSS}
after making the $\alpha$-specialization
by $\rho_{\sigma}$, $M_{\sigma}$, $\bJ_{\sigma}(z)$ and $\bK_{\sigma}$, respectively.
Let $c_1$, $c_2$ and $z_0$ be positive constants depending on $\alpha$ and $\tau$ given in \cite{TracyWidom200?}.
These constants are not explicitly given for $\tau \not= 1$, see Section 1 in \cite{TracyWidom200?}.
Employing the following proposition,
we can prove Theorem \ref{alphashifted} as Theorem \ref{mainthm}
and so we omit the proof.

\begin{prop} \label{cor-sh-limit}
We have
$$
\lim_{n \to \infty} (c_2 n^{\frac{1}{3}})^N \rho_{\sigma} ( \{ [c_1 n +c_2 n^{\frac{1}{3}}x_1], \dots,
[c_1 n + c_2 n^{\frac{1}{3}} x_N] \} ) = \det( K_{\Airy} (x_i,x_j))_{1 \leq i,j \leq N}.
$$
The limit is uniform for $(x_1, \dots, x_N)$ in a compact set of $\bR^N$.
\end{prop}

\begin{proof}
From Theorem \ref{corSS},
we have
$$
\rho_{\sigma} (\{ k_1, \dots, k_N \}) = \sqrt{ \det ( M_{\sigma} (\{k_1, \dots, k_N \}) ) }.
$$
Then we may write the skew matrix $M_{\sigma}$ in the form 
$$
M_{\sigma}= \begin{pmatrix} M_1 & M_2 & \\ - ^tM_2 & M_3 \end{pmatrix},
$$
where we put $N \times N$ matrices $M_1= (\bK_{\sigma}(k_i,k_j))_{1 \leq i,j \leq N}$, 
$M_2=( \bK_{\sigma} (k_i, -k_{N-j+1}))_{1 \leq i,j \leq N}$ 
and $M_3= (\bK_{\sigma}(-k_{N-i+1}, -k_{N-j+1}))_{1 \leq i,j \leq N}$.
Let $D$ be an $N \times N$ diagonal matrix whose $i$-th entry is given by $\bJ_{\sigma}(z_0)^{-1}z_0^i$.
Then $\rho_{\sigma}$ is expressed as
\begin{align*}
\rho_{\sigma} ( \{k_1, \dots, k_N \}) 
&= \sqrt{ (-1)^N \det \begin{pmatrix} M_2 & M_1 \\ M_3 & - ^t M_2 \end{pmatrix} } \\
&= \sqrt{ (-1)^N \det \( \begin{pmatrix} D & 0 \\ 0 & D^{-1} \end{pmatrix} 
\begin{pmatrix} M_2 & M_1 \\ M_3 & - ^t M_2 \end{pmatrix} 
\begin{pmatrix} D^{-1} & 0 \\ 0 & D \end{pmatrix} \) } \\
&= \sqrt{ (-1)^N \det \begin{pmatrix} D M_2 D^{-1} & D M_1 D \\ 
D^{-1} M_3 D^{-1} & - D^{-1}\  ^t M_2 D \end{pmatrix} } \\
&= \Pf \begin{pmatrix} D M_1 D & D M_2 D^{-1} \\ -D^{-1} \ ^t M_2 D & D^{-1} M_3 D^{-1} \end{pmatrix}. 
\end{align*}
Thus we immediately obtain the proposition from the following lemma.
\end{proof} 

\begin{lem} 
We have
\begin{align*}
 \bJ_{\sigma}(z_0)^{-2} z_0^{2c_1 n+c_2 n^{\frac{1}{3}} (x+y)} n^{\frac{1}{3}}
\bK_{\sigma} (c_1 n  + c_2 n^{\frac{1}{3}}x, c_1 n + c_2 n^{\frac{1}{3}} y)  &\to 0, \\ 
z_0^{c_2 n^{\frac{1}{3}} (x-y)} n^{\frac{1}{3}}
\bK_{\sigma} (c_1 n  + c_2 n^{\frac{1}{3}}x, -(c_1 n + c_2 n^{\frac{1}{3}} y)) &\to c_2^{-1} K_{\Airy} (x,y), \\
\bJ_{\sigma}(z_0)^2 z_0^{-(2c_1 n+c_2 n^{\frac{1}{3}} (x+y))} n^{\frac{1}{3}}
\bK_{\sigma} (-(c_1 n  + c_2 n^{\frac{1}{3}}x), -(c_1 n + c_2 n^{\frac{1}{3}} y))    &\to 0, 
\end{align*}
as $n \to \infty$, uniformly in $x$ and $y$ on compact sets in $\bR$.
\end{lem}

The proof of this lemma is obtained by the discussion in Section 6.4 of \cite{TracyWidom200?}.
Since the assumption $\alpha^2 < \tau < \alpha^{-2}$ is not used in that section,
we do not need this assumption in Theorem \ref{alphashifted}.

%
\section{Hall-Littlewood measures}
%

Let $\mP$ be the set of all partitions.
In this section,
we consider the so-called Hall-Littlewood measure on $\mP$, defined by Hall-Littlewood functions.
It is considered as a natural extension of the Schur measure and the shifted Schur measure.
Let $Q_{\lambda}(X;t)$ (respectively $P_{\lambda}(X;t)$) be the Hall-Littlewood $Q$-(respectively $P$-)function
for a partition $\lambda$ (see \cite[III]{Macdonald}).
We define the {\it Hall-Littlewood measure} by
\begin{equation*}
\PP_{\mathrm{HL},X,Y,t} (\{\lambda \}) = \frac1{Z} Q_{\lambda}(X;t) P_{\lambda}(Y;t).
\end{equation*}
Here the constant $Z=Z(X,Y;t)$ is determined by
\begin{equation*}
Z:=\sum_{\lambda \in \mP} Q_{\lambda}(X;t) P_{\lambda}(Y;t) = \prod_{i,j=1}^\infty \frac{1-tX_i Y_j}{1-X_i Y_j},
\end{equation*}
where the second equality is the Cauchy identity for Hall-Littlewood functions.
Since $Q_{\lambda}(X;t)$ is the Schur function $s_{\lambda}(X)$ at $t=0$ 
and the Schur $Q$-function $Q_{\lambda}(X)$ at $t=-1$,
the Hall-Littlewood measure gives the Schur measure at $t=0$ and the shifted Schur measure at $t=-1$.

The mean value and the variance of the size $|\lambda|$ of a partition $\lambda$ are given explicitly as follows.
\begin{thm} \label{MeanValue}
The mean value $\bE (|\lambda|)$ and the variance $\mathrm{Var}(|\lambda|)$ 
of the size $|\lambda|$ of a partition
with respect to the Hall-Littlewood measure
are given by
\begin{align}
\bE(|\lambda|) &= \sum_{k=1}^\infty (1-t^k)p_k(X) p_k(Y), \label{SizeMeanValue} \\
\mathrm{Var}(|\lambda|) &= \sum_{k=1}^\infty k (1-t^k)p_k(X) p_k(Y).
\end{align}
Here $p_k(X)$ is the $k$-th power sum function.
\end{thm}

\begin{proof}
Define a differential operator $\Delta_X$ by 
\begin{equation*}
\Delta_X = \sum_{k=1}^\infty k  p_k(X) \frac{\partial}{\partial p_k(X)}.
\end{equation*}
Since 
$Q_{\lambda}(X;t) 
= \sum_{\rho: |\rho|=|\lambda|} z_{\rho}(t)^{-1} X^{\lambda}_{\rho} (t) \prod_{k=1}^\infty p_k(X)^{m_k(\rho)}$
(see \cite[III-(7.51)]{Macdonald})
we have 
$$
\Delta_X Q_{\lambda}(X;t)= |\lambda| Q_{\lambda}(X;t).
$$
Therefore we obtain
$\bE(|\lambda|) = \frac1{Z} \sum_{\lambda} |\lambda| Q_{\lambda}(X;t) P_{\lambda}(Y;t) = \Delta_X( \log Z)$.
On the other hand,
since $Z= \exp \( \sum_{k=1}^\infty \frac{1-t^k}{k} p_k(X)p_k(Y) \)$ (see \cite[p.223]{Macdonald})
we have $\Delta_X ( \log Z) = \sum_{k=1}^\infty (1-t^k) p_k(X) p_k(Y)$
so that we get \eqref{SizeMeanValue}.

In general, we see that $\bE(|\lambda|^n) = \frac1{Z} \Delta_X^n Z$.
In particular,
it follows that
\begin{align*}
\bE (|\lambda|^2) =& \frac1{Z} \Delta_X^2 (Z) = \frac1{Z} \Delta_X \( Z \frac1{Z} \Delta_X(Z)\)
                  =  \frac1{Z} \Delta_X (Z \Delta_X(\log Z)) \\
                  =&  \frac1{Z} \left\{ \Delta_X (Z) \cdot \Delta_X( \log Z) + Z \Delta_X^2 (\log Z) \right\} \\
                  =& (\Delta_X (\log Z))^2 + \Delta_X^2 (\log Z).
\end{align*}
Therefore we have
\begin{align*}
\mathrm{Var}(|\lambda|) =& \bE (|\lambda|^2) - \bE( |\lambda|)^2
                        =  \Delta_X^2 (\log Z) \\
                        =& \Delta_X \( \sum_{k=1}^\infty (1-t^k) p_k(X) p_k(Y) \) \\
                        =&  \sum_{k=1}^\infty k(1-t^k) p_k(X) p_k(Y). 
\end{align*}
This completes the proof of the theorem.
\end{proof}

\begin{remark}
The mean value $\bE(|\lambda|)$ with respect to the Schur measure is given in \cite{OkounkovUses}.
\qed
\end{remark}

\medskip

Next we consider the mean value $\bE(\lambda_1)$ of $\lambda_1$, the first row of a partition $\lambda$.
Based on the fact in Theorem \ref{MeanValue},
we now examine whether $\bE(\lambda_1)$ has an expression similar to $\bE(|\lambda|)$.
Assume $X=Y$ in the definition of the Hall-Littlewood measure and 
let $M(t,X):= 2\sum_{k=1}^\infty (1-t^k) p_k(X)$.

\begin{example}
The poissonized Plancherel measure for symmetric groups is obtained from the Hall-Littlewood measure 
by putting $t=0$ and the exponential specialization $p_k(X) =p_k(Y)= \sqrt{\xi} \delta_{1k}$.
Hence it follows that $M(t,X)= 2 \sqrt{\xi}$.
Since it is known that $\bE(\lambda_1) \sim 2 \sqrt{\xi}$ as $\xi \to +\infty$ (see e.g. \cite{BOO}),
we have $\lim \frac{\bE(\lambda_1)}{M(t,X)} =1$.
\qed
\end{example}

\begin{example}
The $\alpha$-specialized Schur measure is obtained by putting $t=0$ and the $\alpha$-specialization 
$X=Y= (\overbrace{\alpha, \dots, \alpha}^n, 0, 0, \dots)$.
We have hence $M(t,X)= 2 \sum_{k=1}^\infty n \alpha^k = \frac{2 \alpha}{1-\alpha} n$.
Since it is proved in \cite{Johansson2000} that
$\bE(\lambda_1) \sim  \frac{2 \alpha n}{1-\alpha}$ as $n \to +\infty$,
we have $\lim \frac{\bE(\lambda_1)}{M(t,X)} =1$.
\qed
\end{example}

\begin{example}
The poissonization of the shifted version of the Plancherel measure for symmetric groups is obtained by putting $t=-1$ 
and $p_k(X) =p_k(Y)= \sqrt{\frac{\xi}{2}} \delta_{1k}$.
We have hence $M(t,X) = 4 \sqrt{\frac{\xi}{2}}= 2 \sqrt{2\xi}$.
Since we have proved that $\bE(\lambda_1) \sim 2 \sqrt{2\xi}$ as $\xi \to +\infty$ in Section \ref{Applications},
we have $\lim \frac{\bE(\lambda_1)}{M(t,X)} =1$.
\qed
\end{example}

\begin{example}
The $\alpha$-specialized shifted Schur measure is obtained by putting $t=-1$ and the $\alpha$-specialization 
$X=Y= (\overbrace{\alpha, \dots, \alpha}^n, 0, 0, \dots)$.
Hence we have $M(t, X)= 2 \sum_{k \geq 1: \mathrm{odd}} 2n \alpha^k = \frac{4 \alpha n}{1- \alpha^2}$.
Since it is proved in \cite{TracyWidom200?} that
$\bE(\lambda_1) \sim  \frac{4 \alpha n}{1-\alpha^2}$ as $n \to +\infty$,
we have $\lim \frac{\bE(\lambda_1)}{M(t,X)} =1$.
\qed
\end{example}

\medskip

In view of the examples above,
we might expect that $\frac{\bE (\lambda_1)}{M(t,X)}$ always converges to $1$.
However, we encounter an example that $\frac{\bE (\lambda_1)}{M(t,X)}$ does not converge to $1$ as follows.

\begin{example}
Suppose $0 < t< 1$.
We make the principal specialization $X=Y=(t,t^2, \dots, t^n,0,0,\dots)$ and 
$n \to + \infty$.
Then we have 
$$
Q_{\lambda} = t^{n(\lambda)+|\lambda|}, \qquad 
P_{\lambda}= \frac{t^{n(\lambda)+|\lambda|}}{\prod_{j \geq 1} (t;t)_{m_j(\lambda)}},
$$
where we put $n(\lambda)=\sum_{j \geq 1} (j-1) \lambda_j
= \sum_{j \geq 1} \begin{pmatrix} \lambda_j' \\ 2 \end{pmatrix}$,
$(a;q)_m= \prod_{j=0}^{m-1} (1-aq^j)$
and denote the multiplicity of $j$ in $\lambda$ by $m_j(\lambda)$
(see \cite[III-2 Ex.1]{Macdonald}).
Further we obtain 
\begin{equation*}
Z= \prod_{i,j=1}^\infty \frac{1-t^{i+j+1}}{1-t^{i+j}} = \prod_{r=2}^\infty \frac{1}{1-t^r}.
\end{equation*}
Since 
$2n(\lambda) + |\lambda| = \sum_{j \geq 1} \lambda_j'(\lambda_j'-1) + \sum_{j \geq 1} \lambda_j' 
= \sum_{j \geq 1} (\lambda_j')^2$,
the Hall-Littlewood measure becomes
\begin{equation*}
\PP_{t,\mathrm{Prin}} (\lambda)  
:= \prod_{r=2}^\infty (1-t^r) \frac{t^{\sum_{j \geq 1}(\lambda_j')^2+ |\lambda|}} {\prod_{j \geq 1} (t;t)_{m_j(\lambda)}}.
\end{equation*}
This measure is studied by Fulman \cite{Fulman1999}.

We calculate the distribution function of $\lambda_1$.
Since for a positive integer $h$
$$
 \sum_{\lambda:\lambda_1 <h}  \frac{t^{\sum_{j \geq 1}(\lambda_j')^2+ |\lambda|}} {\prod_{j \geq 1} (t;t)_{m_j(\lambda)}}
= \sum_{\mu:\ell(\mu)<h}  \frac{t^{\sum_{j= 1}^{h-1} (\mu_j)^2 + |\mu|}} {\prod_{j = 1}^{h-1} (t;t)_{\mu_j-\mu_{j+1}}}
= \prod_{\begin{subarray}{c} r=1 \\ r \not\equiv 0, \pm 1 \pmod{2h+1} \end{subarray}}^\infty \frac1{1-t^r},
$$
(the second equality is proved by Andrews \cite{Andrews}, see also \cite{Fulman2000})
we have 
\begin{align}
\PP_{t,\mathrm{Prin}} ( \lambda_1 <h) 
=& \prod_{\begin{subarray}{c} r \geq 2 \\ r \equiv 0, \pm 1 \pmod{2h+1} \end{subarray}} (1-t^r) \notag \\
=& \prod_{k=1}^\infty (1-t^{(2h+1)k})(1-t^{(2h+1)k+1})(1-t^{(2h+1)k-1}). \label{distribution}
\end{align}
The mean value
$\bE (\lambda_1)$ is given by
$\bE (\lambda_1) = \sum_{h=1}^\infty h ( \PP_{t,\mathrm{Prin}}(\lambda_1 < h+1) - \PP_{t,\mathrm{Prin}}(\lambda_1 <h))$.
It follows from \eqref{distribution} that
$\bE (\lambda_1) = t^2+ O(t^3)$ as $t \to +0$.
On the other hand, it is easy to see that
$$
M(t,X) = 2 \sum_{k=1}^\infty (1-t^k) \sum_{j=1}^\infty t^{jk} = 2\sum_{k=1}^\infty (1-t^k) \frac{t^k}{1-t^k}
       = \frac{2t}{1-t}
$$
and therefore $M(t,X) = 2t + 2t^2+ O(t^3)$ as $t \to +0$.
Therefore $\bE(\lambda_1)$ is not equal to $M(t,X)$.
\qed
\end{example}

\medskip

Thus it is interesting to determine when the ratio $\frac{\bE(\lambda_1)}{M(t,X)}$ converges to 1.
We will study this problem in future.

%
%


\noindent
\textsc{Sho Matsumoto}\\
Graduate School of Mathematics, Kyushu University.\\
Hakozaki Fukuoka 812-8581, Japan.\\
e-mail : \texttt{ma203029@math.kyushu-u.ac.jp}\\

\end{document}